\documentclass[letterpaper,10pt,conference]{ieeeconf}
\IEEEoverridecommandlockouts
\usepackage[T1]{fontenc}
\usepackage{amsmath,amssymb,mathtools,graphicx,cite}
\usepackage[hidelinks]{hyperref}
\usepackage{eso-pic}
\newtheorem{assumption}{Assumption}
\newtheorem{theorem}{Theorem}

\newtheorem{remark}{Remark}
\renewenvironment{proof}{\par\noindent\textit{Proof: }}{\hfill\QED\par}
\newcommand{\R}{\mathbb R}
\newcommand{\eps}{\varepsilon}
\newcommand{\diag}{\operatorname{diag}}
\newcommand{\col}{\operatorname{col}}
\newcommand{\dist}{\operatorname{dist}}
\newcommand{\atanh}{\operatorname{atanh}}
\newcommand{\norm}[1]{\left\|#1\right\|}
\graphicspath{{figures/}}

\AddToShipoutPictureFG*{%
 \AtPageUpperLeft{%
  \raisebox{-8mm}[0pt][0pt]{%
   \makebox[\paperwidth][c]{%
    \parbox{\textwidth}{\centering\normalfont\fontsize{7}{8}\selectfont
     This work has been submitted to the IEEE for possible publication. Copyright may be transferred without notice, after which this version may no longer be accessible.}}}}}
\title{\LARGE\bf Prescribed Performance Convex Optimization with Nonlinear Inequality Constraints}
\author{\parbox[t]{0.85\textwidth}{\centering
A.~K.~Gkesoulis, H.~E.~Psillakis, G.~C.~Karras,
C.~P.~Bechlioulis, and K.~P.~Valavanis}%
\thanks{This project is funded by the European Union under Horizon Europe
(grant No. 101136568 - HERON)%
\protect\\[2pt]%
\protect\includegraphics[height=8mm]{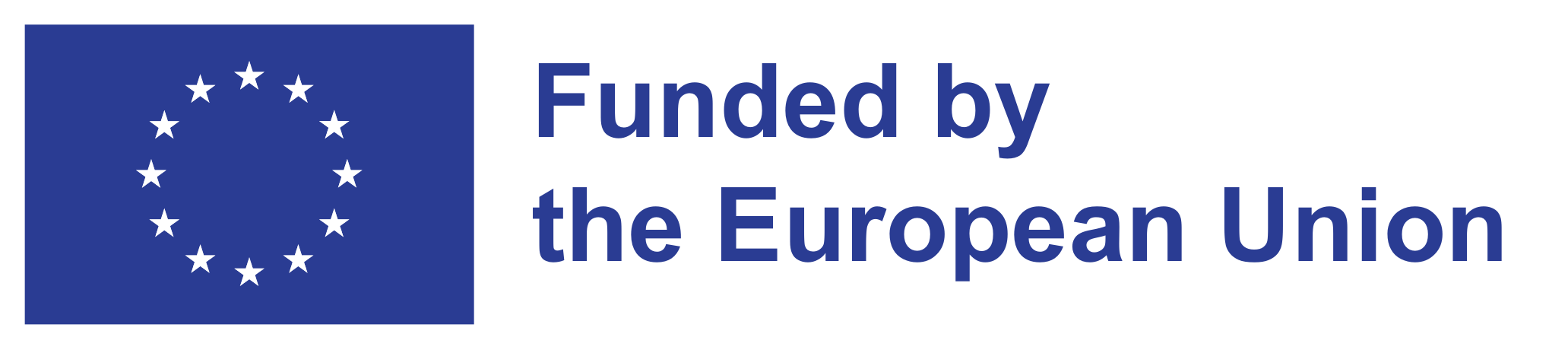}}%
\thanks{A. K. Gkesoulis, G. C. Karras, C. P. Bechlioulis, and K. P. Valavanis
are with HERON - Hellenic Robotics Center of Excellence, Athena Research
Center, Greece (e-mails: \texttt{\{a.gkesoulis, karrasg, chmpechl, k.valavanis\}@athenarc.gr}).}%
\thanks{H. E. Psillakis is with the School of Electrical and Computer
Engineering, National Technical University of Athens, Athens 15780,
Greece (e-mail: \texttt{hpsilakis@central.ntua.gr}).}%
\thanks{G. C. Karras is also with the Department of Informatics and
Telecommunications, University of Thessaly, Lamia 35100, Greece.}%
\thanks{C. P. Bechlioulis is also with the Department of Electrical and
Computer Engineering, University of Patras, Patras 26504, Greece.}%
}
\begin{document}
\maketitle
\thispagestyle{empty}
\pagestyle{empty}

\begin{abstract}
An explicit first-order primal--dual algorithm in continuous time is
developed for strongly convex optimization with smooth convex nonlinear
inequality constraints. From admissible, possibly infeasible initial
states, the algorithm enforces designer-chosen bounds on the stationarity
residual, each constraint violation, and a multiplier-weighted squared
slack that measures complementarity. Under Slater's condition, the
trajectories remain bounded and converge to the exact constrained
optimizer, even when the prescribed bounds have positive limits. The
updates use only function values and first derivatives. When a numerical
lower bound on objective curvature and a strictly feasible point are
known, a constructive design also prescribes an upper bound on the
distance to the optimizer throughout the transient, without knowing the
optimizer or its multipliers. In particular, a chosen accuracy can be
guaranteed by a specified deadline and maintained thereafter.
Simulations illustrate prescribed constraint recovery and optimizer
approach, and compare task completion, entry time, and primal effort
with related methods under equal tuning budgets.
\end{abstract}

\section{Introduction}
Continuous-time projected primal--dual (PD) dynamics solve constrained
optimization through coupled decision and multiplier updated \cite{feijer2010}. Augmented variants admit semi-global
exponential stability for strongly convex objectives and nonlinear convex
inequalities \cite{tang2020}. While these results establish stability and convergence rates, they do not provide separate designer-chosen performance bounds for the optimality errors.

Prescribed performance control assigns transient and steady-state specifications through error transformations
\cite{bechlioulis2014}. Prescribed gradient
or pseudogradient errors for unconstrained optimization and strongly
monotone games are guaranteed in \cite{sun2023} . Work  \cite{garmpis2026} combines prescribed estimator disagreement with gradient tracking and logarithmic  barriers. However, although these approaches prescribe selected gradient or estimator errors, they do not jointly regulate stationarity, constraint satisfaction, 
and complementarity for the original nonlinear constrained problem.

Safe gradient dynamics use velocity quadratic programs (QPs)
\cite{allibhoy2024} or auxiliary multipliers and time-scale separation
\cite{allibhoy2025}. Moving-slack interior-point (IP) methods permit
infeasible initialization through Newton corrections \cite{fazlyab2018}.
Finite- and fixed-time methods \cite{diana2025}, prescribed-time games \cite{qian2024}, and dual-objective speed specifications \cite{tanaka2025} address convergence timing. Direct Karush--Kuhn--Tucker (KKT) residual correction can require
Lagrangian Hessians and saddle-point solves \cite{zenati2026}. Although, these methods address safety or convergence timing, they do not combine separate prescribed performance bounds for the
these errors with an explicit first-order primal--dual update.

Motivated by the discussion above, a first-order primal--dual algorithm
is developed for prescribed performance convex optimization with nonlinear
inequality constraints. The algorithm enforces prescribed bounds on the
optimality errors, including constraint violations, and converges
asymptotically to the exact constrained optimizer. The contributions can be summarized as:\begin{itemize}
\item An explicit first-order algorithm is developed to enforce separate
prescribed bounds on the stationarity residual, constraint violation,
and multiplier-weighted squared slack. It converges to the exact
constrained optimizer even from admissible infeasible starts.
\item A strong-convexity certificate \cite{mangasarian1988} is utilized to
guarantee a designer-chosen prescribed distance bound and accuracy deadline in the case of 
known curvature and Slater data, without knowing the optimizer or
its multipliers.
\item A common-task comparison is provided under equal tuning budgets.
At the intermediate speed cap, the proposed method completes more
tasks than PD and safe gradient. On the same six cases completed by
IP, it achieves lower mean entry time and primal effort.
\end{itemize}
The rest of the paper is organized as follows.
Sections~\ref{sec:problem}--\ref{sec:design} define the problem and
algorithm, Sections~\ref{sec:analysis}--\ref{sec:primal} establish the
guarantees, and Section~\ref{sec:numerical} reports simulations.

\noindent\textit{Notation:}
The sets $\R$ and $\mathbb N$ denote real numbers and positive integers.
For $u\in\R$, $u_+:=\max\{u,0\}$,  vector positive parts are
componentwise. Norms are Euclidean. The symbols $\col$, $\diag$,
$I_k$ ($k\in\mathbb N$), and $\odot$ denote column stacking, diagonal
construction, the identity matrix, and componentwise multiplication.
For a point $y$ and a nonempty set $A$ in the same Euclidean space,
$\dist(y,A):=\inf_{a\in A}\norm{y-a}$.

\section{Problem formulation}\label{sec:problem}
Let $n,m\in\mathbb N$.
For $f:\R^n\to\R$ and $g_i:\R^n\to\R$, consider
\begin{equation}\label{eq:problem}
 \begin{aligned}
 \underset{x\in\R^n}{\operatorname{minimize}}\quad & f(x)\\
 \text{subject to}\quad &g_i(x)\leq0,\quad i=1,\ldots,m .
 \end{aligned}
\end{equation}
The analysis uses the following assumptions.
\begin{assumption}\label{ass:problem}
The functions $f,g_i$ are $C^2$.
For some $\mu>0$, $\nabla^2f(x)\succeq\mu I_n$ for all $x\in\R^n$;
each $g_i$ is convex.
\end{assumption}
\begin{assumption}\label{ass:slater}
A Slater point $x_s\in\R^n$ exists: $g_i(x_s)<0$ for all $i$.
\end{assumption}
Slater's condition ensures that the feasible set is nonempty. It is
closed and convex because the constraint functions are continuous and
convex. Strong convexity makes $f$ coercive and strictly convex, so the
problem admits a unique optimizer $x^\star$. Slater's condition ensures the existence of optimal
multipliers by convex duality \cite{feijer2010}.
Define $g(x):=\col(g_1(x),\ldots,g_m(x))$. Its Jacobian
$J_g:\R^n\to\R^{m\times n}$ has rows $\nabla g_i(x)^\top$.
For $\lambda\in\R_{\geq0}^m$, let
$\mathcal L(x,\lambda):=f(x)+\lambda^\top g(x)$ and define
\begin{equation}\label{eq:residual}
 r(x,\lambda):=\nabla f(x)+J_g(x)^\top\lambda .
\end{equation}
The KKT set $\mathcal K:=\{x^\star\}\times\mathcal M^\star$ has
$\mathcal M^\star:=\{\lambda\geq0:r(x^\star,\lambda)=0,\
\lambda_i g_i(x^\star)=0\ \forall i\}$.

\section{Prescribed performance specifications}\label{sec:specifications}
Along $x(t)\in\R^n$ and $\lambda(t)\in\R_{\geq0}^m$, write
$r(t):=r(x(t),\lambda(t))$ for the stationarity residual and
$g_i^+(x):=(g_i(x))_+$ for the constraint violation, with
$g^+(x):=(g(x))_+$. Define $s_i(t):=(-g_i(x(t)))_+$,
$s(t):=\col(s_i(t))$, and $C_-(t):=\sum_i\lambda_i(t)s_i(t)$.
The multiplier-weighted squared slack is defined as
\begin{equation}\label{eq:complementarity}
 S(t):=\sum_{i=1}^m\lambda_i(t)s_i(t)^2 .
\end{equation}
Only strictly satisfied constraints contribute to $S$.

Let $\mathcal P$ be the class of prescribed performance functions
$\rho\in C^1([0,\infty),(0,\infty))$ that are nonincreasing,
have a positive limit $\rho_\infty:=\lim_{t\to\infty}\rho(t)>0$, and satisfy
$\sup_{t\geq0}|\dot\rho(t)|<\infty$.
Consequently $\rho_\infty\leq\rho(t)\leq\rho(0)$ and
$\int_0^\infty|\dot\rho(t)|dt=\rho(0)-\rho_\infty$.
Choose $\rho_{1,j},\rho_{2,i},\rho_S\in\mathcal P$ for
$j=1,\ldots,n$ and $i=1,\ldots,m$. The goal of this paper is to provide a dual-primal design that achieves prescribed performance convex optimization with nonlinear inequality constraints, i.e. to bound the quantities 
$|r_j(t)|$, $g_i^+(x(t))$, and $S(t)$ by the prescribed performance functions defined beforehand, together with $x(t)\to x^\star$.
To ensure the problem is well defined, the initial values of the performance functions must satisfy
\begin{equation}\label{eq:initial}
 \begin{aligned}
 |r_j(0)|<\rho_{1,j}(0),~
 g_i^+(x_0)<\rho_{2,i}(0),~
 S(0)<\rho_S(0),
 \end{aligned}
\end{equation}
where $x_0:=x(0)$ and $\lambda_0:=\lambda(0)\geq0$.  

\section{Primal--dual algorithm}\label{sec:design}
For the error transformations, set
$\rho_1(t):=\col(\rho_{1,j}(t))$,
$\rho_2(t):=\col(\rho_{2,i}(t))$,
$R_1(t):=\diag(\rho_{1,j}(t))$, and
$R_2(t):=\diag(\rho_{2,i}(t))$. Define
$\xi_1(t):=R_1(t)^{-1}r(t)$ and
$\xi_2(t):=R_2(t)^{-1}g(x(t))$.
On $|\xi_{1,j}(t)|<1$, set
\begin{equation}\label{eq:eps1}
 \begin{aligned}
 \eps_{1,j}(t)&:=\atanh(\xi_{1,j}(t)),\\
 D_{1,jj}(t)&:=(1-\xi_{1,j}(t)^2)^{-1},\quad j=1,\ldots,n,
 \end{aligned}
\end{equation}
where $\atanh u:=\tfrac12\log((1+u)/(1-u))$ for $u\in(-1,1)$.
For $\xi_{2,i}(t)<1$, define
\begin{equation}\label{eq:eps2}
 \begin{aligned}
 \eps_{2,i}(t)&:=\begin{cases}0,&\xi_{2,i}(t)\leq0,\\
 -\log(1-\xi_{2,i}(t)),&\xi_{2,i}(t)>0,\end{cases}\\
 D_{2,ii}(t)&:=\begin{cases}0,&\xi_{2,i}(t)\leq0,\\
 (1-\xi_{2,i}(t))^{-1},&\xi_{2,i}(t)>0,\end{cases}
 \quad i=1,\ldots,m.
 \end{aligned}
\end{equation}
Stack $\eps_1:=\col(\eps_{1,1},\ldots,\eps_{1,n})$ and
$\eps_2:=\col(\eps_{2,1},\ldots,\eps_{2,m})$, and set
$D_1:=\diag(D_{1,11},\ldots,D_{1,nn})$ and
$D_2:=\diag(D_{2,11},\ldots,D_{2,mm})$.
The choice in \eqref{eq:eps2} keeps $D_2\eps_2$ locally Lipschitz
in $\xi_2$ even though $D_2$ jumps at zero.
For $0\leq S(t)<\rho_S(t)$, set
\begin{equation}\label{eq:epsc}
 \begin{aligned}
 \eps_c(t)&:=-\log(1-S(t)/\rho_S(t)),\\
 \chi(t)&:=\eps_c(t)/(\rho_S(t)-S(t)).
 \end{aligned}
\end{equation}
The proposed algorithm is
\begin{equation}\label{eq:controller}
 {\begin{aligned}
 \dot x(t)&=-k_1D_1(t)R_1(t)^{-1}\eps_1(t),\\
 \dot\lambda_i(t)&=k_\lambda(R_2(t)^{-1}D_2(t)\eps_2(t))_i-2k_1\chi(t)\lambda_i(t)s_i(t),
 \end{aligned}}
\end{equation}
where $k_1,k_\lambda$ are positive gains. The first multiplier term corrects positive constraint error, while the
second damps multipliers of strictly satisfied constraints.

\section{Performance Analysis}\label{sec:analysis}
The next theorem gives global existence, prescribed error bounds, and convergence to the optimizer.
\begin{theorem}\label{thm:main}
Under Assumptions~\ref{ass:problem}--\ref{ass:slater}, choose performance
functions in $\mathcal P$ satisfying \eqref{eq:initial} with
$\lambda(0)\geq0$. Then \eqref{eq:controller} has a unique global
classical solution with $\lambda(t)\geq0$ for all $t\geq0$.
The following hold:
\begin{enumerate}
\item For every $t\geq0$, $|r_j(t)|<\rho_{1,j}(t)$,
$g_i^+(x(t))<\rho_{2,i}(t)$, and $S(t)<\rho_S(t)$.
The quantities $x$, $\lambda$, $\eps_1$, $\eps_2$, $\eps_c$, $\chi$,
$\dot x$, $\dot\lambda$, $D_1R_1^{-1}\eps_1$, and
$R_2^{-1}D_2\eps_2$ are bounded.
\item As $t\to\infty$, $r$, $g^+$, $C_-$, $S$, $\eps_1$, $\eps_2$,
$\eps_c$, $\dot x$, and $\dot\lambda$ tend to zero. Moreover,
\begin{equation}\label{eq:limits}
 \lim_{t\to\infty}x(t)=x^\star,~
 \lim_{t\to\infty}\dist((x(t),\lambda(t)),\mathcal K)=0.
\end{equation}

If also $\mathcal M^\star=\{\lambda^\star\}$, then
$\lim_{t\to\infty}\lambda(t)=\lambda^\star$.
\end{enumerate}
\end{theorem}
\begin{proof}
\emph{Local solution and state bounds.}
Allow $\lambda\in\R^m$ on the open state domain
$|r_j|<\rho_{1,j}$, $g_i<\rho_{2,i}$, $S<\rho_S$.
The formula for $\eps_c$ also applies to $S<0$.
Local Lipschitz continuity gives a unique maximal solution on
$[0,t_{\max})$, where these domain inequalities hold. Along that
solution each multiplier equation has a continuous coefficient of
$\lambda_i$ and a nonnegative input. Its integrating factor gives
$\lambda_i\geq0$, hence $S\geq0$.

Choose a Slater point $x_s$ and set $\delta_i:=-g_i(x_s)$,
$\delta:=\min_{1\leq i\leq m}\delta_i>0$ and
$G:=\norm{\nabla f(x_s)}$. In $r^\top(x-x_s)$, the objective term is
bounded below by $\mu\norm{x-x_s}^2-G\norm{x-x_s}$ by strong convexity.
Each constraint gradient satisfies
$\nabla g_i(x)^\top(x-x_s)\geq g_i(x)+\delta_i\geq\delta-s_i$.
Using $\lambda_i\geq0$ and Cauchy--Schwarz on $r^\top(x-x_s)$ gives
\begin{equation}\label{eq:slater-inequality}
 \mu\norm{x-x_s}^2+\delta\sum_i\lambda_i
 \leq(\norm r+G)\norm{x-x_s}+C_- .
\end{equation}
By weighted Cauchy--Schwarz and Young's inequality,
$C_-\leq\sqrt{S\sum_i\lambda_i}\leq
\tfrac{\delta}{2}\sum_i\lambda_i+S/(2\delta)$.
With $\norm r\leq\norm{\rho_1(0)}$ and $S<\rho_S(0)$,
square completion in \eqref{eq:slater-inequality} now bounds
$x$ and $\sum_i\lambda_i$, hence also $C_-$. Jensen's inequality for
$u^{3/2}$ with weights $\lambda_i/\sum_j\lambda_j$ gives
$\sum_i\lambda_i s_i^3\geq c S^{3/2}$ for some $c>0$ by boundedness
of the multiplier sum. If the sum is zero, both sides vanish.

\emph{Lyapunov derivative.}
Set
\begin{equation}\label{eq:V}
 V:=\tfrac12\norm{\eps_1}^2+
     \tfrac{k_\lambda}{2k_1}\norm{\eps_2}^2+\tfrac12\eps_c^2.
\end{equation}
The derivative of $\eps_{2,i}^2/2$ with respect to $\xi_{2,i}$
is $(D_2\eps_2)_i$, including at zero where both one-sided
derivatives vanish. Using \eqref{eq:controller}, the chain rule gives
\begin{align*}
 \frac{d}{dt}\frac{\norm{\eps_1}^2}{2}
  &=-\frac{\dot x^\top}{k_1}(\dot r-\dot R_1\xi_1),\\
 \frac{d}{dt}\frac{k_\lambda\norm{\eps_2}^2}{2k_1}
  &=\frac{k_\lambda}{k_1}\eps_2^\top D_2R_2^{-1}
     (J_g\dot x-\dot R_2\xi_2),\\
 \frac{d}{dt}\frac{\eps_c^2}{2}
  &=\chi\bigl(\dot S-\tfrac{S}{\rho_S}\dot\rho_S\bigr).
\end{align*}
Here $\dot r=H_{\mathcal L}\dot x+J_g^\top\dot\lambda$ with
$H_{\mathcal L}:=\nabla^2 f+\sum_i\lambda_i\nabla^2g_i\succeq\mu I_n$.
Differentiating $S$ uses $s_i^2(R_2^{-1}D_2\eps_2)_i=0$ and
$\frac{d}{dt}s_i^2=-2s_i\nabla g_i^\top\dot x$, giving
\begin{equation}\label{eq:storage-derivative}
 \dot S=-2k_1\chi\sum_i\lambda_i s_i^3
        -2\dot x^\top J_g^\top(\lambda\odot s).
\end{equation}
The correction and damping contributions in
$-\dot x^\top J_g^\top\dot\lambda/k_1$ cancel, respectively,
the $J_g\dot x$ term in the violation energy and the
$\dot x$ term in $\chi\dot S$. Summation therefore yields
\begin{equation}\label{eq:Vdot}
 \begin{aligned}
 \dot V={}&-\tfrac1{k_1}\dot x^\top H_{\mathcal L}\dot x
   -2k_1\chi^2\sum_i\lambda_i s_i^3+\tfrac1{k_1}\dot x^\top\dot R_1\xi_1\\
 &
 -\tfrac{k_\lambda}{k_1}\eps_2^\top D_2R_2^{-1}\dot R_2\xi_2-\chi\tfrac{S}{\rho_S}\dot\rho_S .
 \end{aligned}
\end{equation}

\emph{Boundedness and continuation.}
Where $(D_2\eps_2)_i>0$, $0<\xi_{2,i}<1$. Thus the
$\dot R_2$ term in \eqref{eq:Vdot} is bounded above by
$\tfrac{k_\lambda}{k_1}\max_i\sup_{t\geq0}(-\dot\rho_{2,i}(t))$
times $\sum_i(R_2^{-1}D_2\eps_2)_i$. Choose
$a_F>k_1^{-1}\max_i\sup_{t\geq0}(-\dot\rho_{2,i}(t))$
and set $F:=V-a_F\sum_i\lambda_i$. The correction in $\dot F$ leaves a negative multiple of the same
sum and adds $2k_1a_F\chi C_-$.

Young's inequality absorbs the bounded-rate $\dot R_1$ term into half
the primal dissipation plus a constant. Bounded $C_-,\dot\rho_S$ and
$S/\rho_S<1$ bound the storage terms by
$-2k_1\chi^2\sum_i\lambda_i s_i^3+C\chi$.
If $S/\rho_S\geq1/2$, the positive limit of $\rho_S$ bounds the cubic
sum away from zero. Since $\chi\geq\eps_c/\rho_S(0)$, completing the
square gives an upper bound of the form $-c\eps_c^2+C$.
Otherwise $\eps_c\leq\log2$ and $\chi$ is bounded by
$2\log2/\inf_{t\geq0}\rho_S(t)$, therefore this form  holds.

The remaining negative terms control the other energies:
$D_1\succeq I_n$ and bounded $R_1$ make $\norm{\dot x}^2$
dominate the positive multiple of $\norm{\eps_1}^2$. 
$(D_2\eps_2)_i=\eps_{2,i}e^{\eps_{2,i}}\geq\eps_{2,i}^2$
and bounded $R_2$ give the corresponding bound for
$\sum_i(R_2^{-1}D_2\eps_2)_i$. Thus, for constants $c>0$ and $C<\infty$,
$\dot F\leq-cV+C\leq-cF+C$. Scalar comparison bounds $F$ above,
and $V-F=a_F\sum_i\lambda_i$ is bounded, so $V$ is bounded.
Inverting the transforms gives uniform margins from the singular
boundaries. These margins and the positive profile limits bound
both update products and $\chi$, and \eqref{eq:controller} bounds
the velocities. At a finite $t_{\max}$ the state would have an
interior limit and admit a local extension. Thus the solution is global.

\emph{Exact convergence.}
The profile derivatives in \eqref{eq:Vdot} are integrable by
finite total variation and their coefficients are bounded.
Since $V\geq0$, integration gives
$ \int_0^\infty\!
 \left(\norm{\dot x}^2+\chi^2\sum_i\lambda_i s_i^3\right)dt<\infty .$
Bounded derivatives make the states and profiles Lipschitz in time.
The integrands are continuous functions of these quantities on a compact
set away from the singular boundaries, hence uniformly continuous in time. Since
$\chi\geq S/\rho_S(0)^2$, the cubic sum estimate bounds
$\chi^2\sum_i\lambda_i s_i^3$ below by a positive multiple of
$S^{7/2}$. Barbalat's lemma gives $\dot x\to0$ and $S\to0$,
so $\eps_1,r\to0$ by the primal update and $C_-\to0$ by the
weighted estimate.

Bounded constraint gradients now give $\frac{d}{dt}g_i(x(t))\to0$.
If $g_i(x(t))\geq2\eta$ for some $\eta>0$ at arbitrarily late
times, it remains above $\eta$ on a following interval of any
fixed length for sufficiently late such times. There $s_i=0$ and
$\dot\lambda_i\geq k_\lambda\eta/\rho_{2,i}(0)^2$ by
$-\log(1-u)\geq u$ for $0<u<1$. Taking the interval sufficiently
long contradicts the multiplier bound. Thus $g^+\to0$.
Since $C_-\to0$ and the multipliers are bounded,
$\lambda_i g_i(x)\to0$. Each cluster pair is therefore KKT.
Boundedness and uniqueness of the primal optimizer give
\eqref{eq:limits}, and a unique optimal multiplier gives
$\lambda\to\lambda^\star$. Finally the positive profile limits
give $\eps_2,\eps_c,\chi\to0$, so $\dot\lambda\to0$.
\end{proof}
\begin{remark}
The regulated errors converge to zero even for constant performance
functions. 
\end{remark}

\section{Prescribed approach to the optimizer}\label{sec:primal}
The next result prescribes a bound on optimizer distance.    
\begin{theorem}\label{thm:primal}
Under the hypotheses of Theorem~\ref{thm:main}, suppose that a numerical
lower bound $\mu>0$ on the curvature of $f$ and a Slater point $x_s$
are known. Set $\delta_i:=-g_i(x_s)>0$,
$\delta:=\min_{1\leq i\leq m}\delta_i$,
$G:=\norm{\nabla f(x_s)}$, and $K_s:=G^2/(4\mu)$.
For admissible performance functions, let $b_\rho(t):=\norm{\rho_1(t)}$
and $A_\rho(t):=(b_\rho(t)+G)^2/(4\mu)$, and define
\begin{equation*}
 \begin{aligned}
 \bar C_\rho(t)&:=\frac{\rho_S(t)+
 \sqrt{\rho_S(t)^2+4\delta A_\rho(t)\rho_S(t)}}{2\delta},\\
 B_\rho(t)&:=\bar C_\rho(t)+K_s
       \max_{1\leq i\leq m}\frac{\rho_{2,i}(t)}{\delta_i},\\
 \bar d_\rho(t)&:=\frac{b_\rho(t)+
      \sqrt{b_\rho(t)^2+4\mu B_\rho(t)}}{2\mu}.
 \end{aligned}
\end{equation*}
Then $\norm{x(t)-x^\star}<\bar d_\rho(t), $  $ \forall t\geq0$, and $\lim_{t\to\infty}x(t)= x^\star.$
\end{theorem}
\begin{proof}
Fix $\lambda^\star\in\mathcal M^\star$. Strong convexity bounds
$\mu\norm{x-x^\star}^2$ by
$(\nabla f(x)-\nabla f(x^\star))^\top(x-x^\star)$.
Substitute $\nabla f(x)=r(x,\lambda)-J_g(x)^\top\lambda$ and
$\nabla f(x^\star)=-J_g(x^\star)^\top\lambda^\star$.
By convexity,
$J_g(x)(x-x^\star)\geq g(x)-g(x^\star)$ and
$J_g(x^\star)(x-x^\star)\leq g(x)-g(x^\star)$, componentwise.
Using $\lambda,\lambda^\star\geq0$,
$\lambda^\top g(x^\star)\leq0$ and
$(\lambda^\star)^\top g(x^\star)=0$ yields the signed estimate
\begin{equation}\label{eq:certificate}
 \begin{aligned}
 \mu\norm{x-x^\star}^2
 \leq r(x,\lambda)^\top(x-x^\star)-\lambda^\top g(x)+(\lambda^\star)^\top g(x).
 \end{aligned}
\end{equation}
Write $\delta_i:=-g_i(x_s)>0$, $\delta:=\min_i\delta_i$, and
$G:=\norm{\nabla f(x_s)}$. At $(x_s,0)$, \eqref{eq:certificate} gives
$\mu\norm{x_s-x^\star}^2+\sum_i\delta_i\lambda_i^\star
\leq G\norm{x_s-x^\star}$, so
$\sum_i\delta_i\lambda_i^\star\leq G^2/(4\mu)$ by completing the square.
Thus along the solution,
$(\lambda^\star)^\top g(x)\leq(\lambda^\star)^\top g^+(x)
\leq (G^2/(4\mu))\max_i(\rho_{2,i}/\delta_i)$.

The pointwise estimate \eqref{eq:slater-inequality} and square completion
give $\delta\sum_i\lambda_i\leq(\norm{\rho_1}+G)^2/(4\mu)+C_-$.
By weighted Cauchy--Schwarz, $C_-^2\leq S\sum_i\lambda_i$; since
$\norm{\rho_1}>0$ and $S<\rho_S$, it follows that
\[
 \delta C_-^2<\rho_S\left(\frac{(\norm{\rho_1}+G)^2}{4\mu}+C_-\right).
\]
Solving this inequality bounds $C_-$ strictly by the positive root.
Using $-\lambda^\top g(x)\leq C_-$ and $\norm r\leq\norm{\rho_1}$
in \eqref{eq:certificate} now gives, with time suppressed,
\begin{equation}\label{eq:general-distance-bound}
 \begin{aligned}
 \mu\norm{x-x^\star}^2
 &< \norm{\rho_1}\norm{x-x^\star}
    +\frac{G^2}{4\mu}\max_i\frac{\rho_{2,i}}{\delta_i}\\
 &\quad+\frac{\rho_S+\sqrt{\rho_S^2+
   \delta\rho_S(\norm{\rho_1}+G)^2/\mu}}{2\delta}.
 \end{aligned}
\end{equation}
Solving the inequality 
for $\norm{x(t)-x^\star}$ gives the bound $\bar d_\rho$ and
exact
convergence follows from Theorem~\ref{thm:main}.
\end{proof}
\begin{remark}\label{rem:allocation}
For a designer-chosen $\rho_x\in\mathcal P$, an admissible allocation
with $\bar d_\rho(t)\leq\rho_x(t)$ for every $t\geq0$
 guarantees
\begin{equation}\label{eq:primal-tube}
{\norm{x(t)-x^\star}<\rho_x(t),\quad t\geq0.}
\end{equation}
Under this allocation, if $\rho_x(T)\leq\varepsilon$ for a chosen
deadline $T>0$ and tolerance $\varepsilon>0$, then
$\norm{x(t)-x^\star}<\varepsilon$ for every $t\geq T$. This allocation realizes a prescribed bound on optimizer distance.
For $G>0$, take $\theta_r,\theta_c,\theta_g>0$ summing to one.
With time suppressed, choose
\begin{equation}\label{eq:explicit-allocation}
 \begin{aligned}
 \rho_{1,j}&:=\theta_r\mu\rho_x/\sqrt n,\\
 \rho_{2,i}&:=4\theta_g\mu^2\delta_i\rho_x^2/G^2,\\
 \rho_S&:=\frac{\delta\theta_c^2\mu^2\rho_x^2}
 {(\theta_r\mu+G/\rho_x)^2/(4\mu)+\theta_c\mu}.
 \end{aligned}
\end{equation}
Substitution gives $\theta_c\mu\rho_x^2$ as the positive root of the
quadratic in $C_-$. In \eqref{eq:general-distance-bound}, the violation
term is $\theta_g\mu\rho_x^2$ and the distance coefficient is
$\theta_r\mu\rho_x$.
Substituting $\rho_x$ for the distance gives equality, so
$\bar d_\rho=\rho_x$. Replacing $\delta_i$ by $\delta$ in $\rho_{2,i}$
gives the same bound, since $\max_i\delta/\delta_i=1$.
If $G=0$, the multiplier bound gives $\lambda^\star=0$.
Use $\theta_r,\theta_c>0$ summing to one and the same $\rho_1,\rho_S$
formulas, with any $\rho_{2,i}\in\mathcal P$ satisfying \eqref{eq:initial}.
As a function of $\rho_x>0$, the numerator of $\rho_S$ increases
and its positive denominator is nonincreasing. This smooth function
has bounded derivative on the compact positive range of $\rho_x$,
therefore the allocated profiles are in $\mathcal P$.
Check \eqref{eq:initial} at $\rho_x(0)$.
The distance need not decrease monotonically. Separate estimates of
negative complementarity and weighted violation can make the bound
conservative. The allocation needs neither $x^\star$ nor $\lambda^\star$.
\end{remark}

\begin{remark}\label{rem:objective-value}
The same allocation bounds the objective-value error on both sides.
Optimality of $x^\star$ for $\mathcal L(\cdot,\lambda^\star)$ gives
$f(x)-f(x^\star)\geq-(\lambda^\star)^\top g^+(x)$.
For the current $\lambda\geq0$, strong convexity and square completion
give $\mathcal L(x,\lambda)-\mathcal L(x^\star,\lambda)
\leq\norm r^2/(2\mu)$. Since $\mathcal L(x^\star,\lambda)\leq f(x^\star)$
and $-\lambda^\top g(x)\leq C_-$, the upper bound is
$f(x)-f(x^\star)\leq\norm r^2/(2\mu)+C_-$.
The allocation in Remark~\ref{rem:allocation} therefore yields
\begin{equation}
 \begin{aligned}
 -\theta_g\mu\rho_x(t)^2
 \leq f(x(t))-f(x^\star)\leq \bigl(\theta_r^2/2+\theta_c\bigr)\mu\rho_x(t)^2 .
 \end{aligned}
\end{equation}
Use $\theta_g=0$ here if $G=0$. If $\rho_x(T)\leq\varepsilon$,
both bounds hold with $\rho_x(t)^2$ replaced by $\varepsilon^2$
for every $t\geq T$.
\end{remark}

\section{Numerical results}\label{sec:numerical}
The simulations illustrate prescribed constraint satisfaction and optimizer
approach, and compare the proposed algorithm with PD \cite{feijer2010},
safe gradient \cite{allibhoy2024}, and IP \cite{fazlyab2018}.
\subsection{Prescribed constraint satisfaction}
\begin{figure*}[t]
 \centering
 \includegraphics[width=0.85\textwidth]{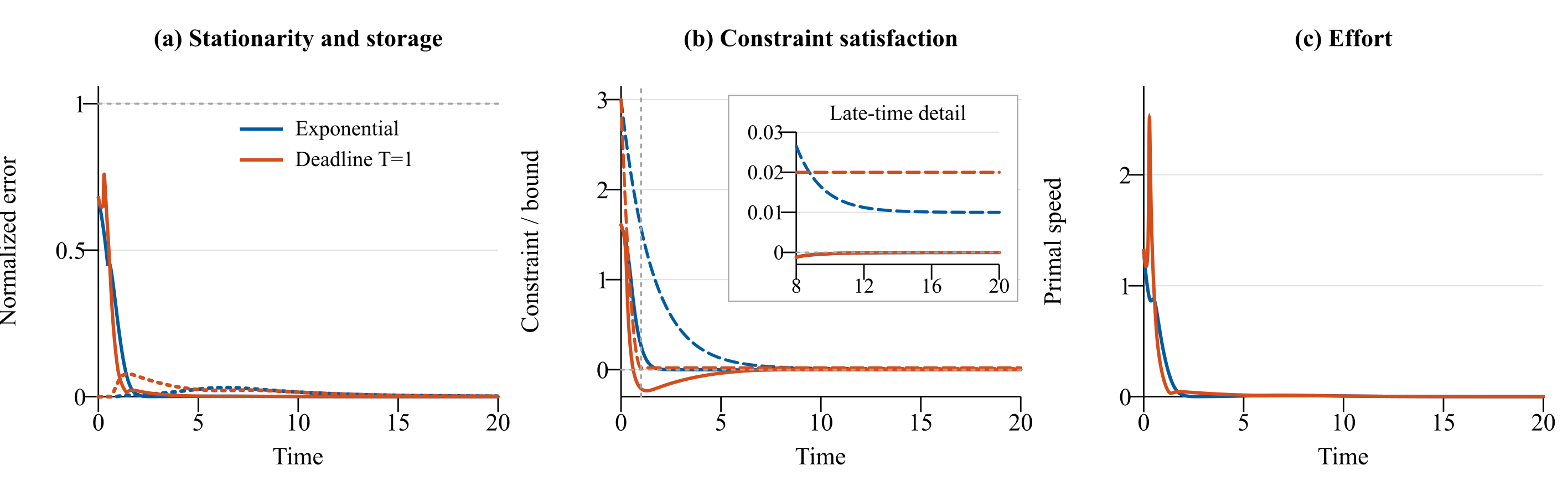}
 \caption{Exponential and one-second choices of $\rho_{2,1}$.
 Left: $\max_j |r_j|/\rho_{1,j}$ (solid) and $S/\rho_S$ (dotted).
 Middle: $g_1$ (solid) and $\rho_{2,1}$ (dashed), with a late-time inset.
 Right: primal speed.}
 \label{fig:specifications}
\end{figure*}

This example isolates the effect of the performance function for the
first constraint. The proposed algorithm is run twice from the same
infeasible initial state, with identical gains and all other performance
functions fixed. The two runs compare exponential recovery with a
prescribed violation tolerance to be met by $t=1$. For $x:=\col(x_1,x_2)\in\R^2$, consider
\begin{equation}\label{eq:benchmark}
 \begin{aligned}
 f(x)&:=\tfrac12((x_1-2)^2+x_2^2)
       +\tfrac14((x_1-1)^4+x_2^4)\\
 &\quad+\tfrac1{40}(x_1+x_2-1)^4,\\
 g_1(x)&:=x_1^2+x_2^2-1,\\
 g_2(x)&:=(x_1-1.1)^2+(x_2-.8)^2-.83^2 .
 \end{aligned}
\end{equation}
At $\col(1,0)$, $g=\col(0,-.0389)$ and $\nabla f=-.5\nabla g_1$.
Thus $\lambda^\star=\col(.5,0)$ verifies the KKT conditions, giving
$x^\star=\col(1,0)$. Use $x_0:=\col(1.5,-.6)$, $\lambda_0:=0$,
and $k_1=k_\lambda:=1$. No Slater or curvature data
are used to select
$\rho_{1,1}(t):=.06+1.14e^{-.55t}$,
$\rho_{1,2}(t):=.08+1.12e^{-.45t}$,
$\rho_{2,2}(t):=.02+2.98e^{-.35t}$, and
$\rho_S(t):=.001+.999e^{-.55t}$.
The runs differ only in their first-constraint specifications $\rho_{2,1}^{\rm E}(t):=.01+2.99e^{-.65t}$, 
 $\rho_{2,1}^{\rm D}(t):=.02+2.98(1-t)_+^2$. Theorem~\ref{thm:main} gives $g_1^+(x(t))<.02$ for all $t\geq1$
with the second choice. In Fig.~\ref{fig:specifications}, we can observe that both runs
respect their bounds. The one-second specification recovers faster
from the violation and increases peak speed from $1.320$ to $2.518$.

\subsection{Prescribed optimizer-distance performance}
This example illustrates how to prescribe an accuracy deadline for the
distance to the optimizer while keeping the benchmark, initial state,
and gains fixed. For each deadline $T\in\{2,4,8\}$, the distance profile
is chosen as $\rho_x(t):=.02+3.98(.03/3.98)^{t/T}$, 
so that $\rho_x(0)=4$ and $\rho_x(T)=.05$.

\begin{figure*}[t]
 \centering
 \includegraphics[width=0.85\textwidth]{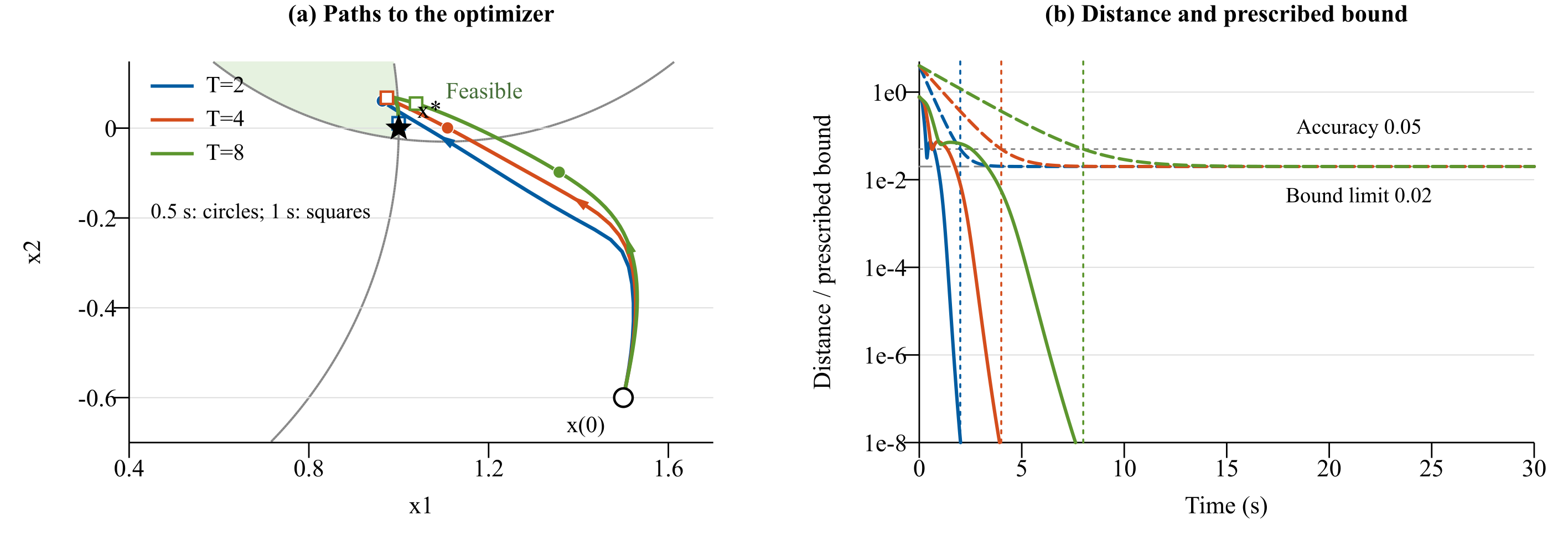}
 \caption{Optimizer approach at fixed gains. Left: paths and time markers,
 feasible region (shaded), and optimizer (star). Right: distances
 (solid), bounds (dashed), and deadlines through $t=30$.
 Distances below $10^{-8}$ are outside the plotted range.}
 \label{fig:primal}
\end{figure*}

The internal performance functions are assigned by
Remark~\ref{rem:allocation} using the certified curvature bound $\mu:=1$
and the Slater point $x_s:=\col(.8,.3)$, which gives $\delta_1=.27$,
$\delta_2=.3489$, and $G^2=1.56601682$. The weights are
$(\theta_r,\theta_c,\theta_g):=(.4,.3,.3)$, with the common margin
$\delta=.27$ used in $\rho_{2,i}$. The resulting functions satisfy
\eqref{eq:initial} and give $\bar d_\rho(t)=\rho_x(t)$. By
Theorem~\ref{thm:primal}, the distance $\norm{x(t)-x^\star}$ stays
below $\rho_x(t)$ for all $t\geq0$ and below $.05$ for all $t\geq T$.
Figure~\ref{fig:primal} shows the trajectories and their distance
profiles. In these runs, shorter prescribed deadlines led to earlier
entry into the same accuracy neighborhood.

\subsection{Comparative study}
This study compares the proposed algorithm with projected
PD~\cite{feijer2010}, safe gradient~\cite{allibhoy2024}, and
prediction--correction IP~\cite{fazlyab2018}. The aim is to test which
methods reach a common optimizer tolerance while respecting the same
constraint violation and speed limits. Parameters are selected on three
development problems and evaluated without retuning on six test problems
generated with disjoint seeds. The test problems have dimensions
$2,5,10$, two active quadratic constraints at the optimizer, and strongly
convex objectives with a quartic term $.005\sum_j x_j^4$ and a quadratic
matrix of condition number $3$ or $10$. All methods share $x_0$ on each
problem and PD~\cite{feijer2010} and the proposed method start with zero
multipliers.

\begin{figure*}[t]
 \centering
 \includegraphics[width=0.85\textwidth]{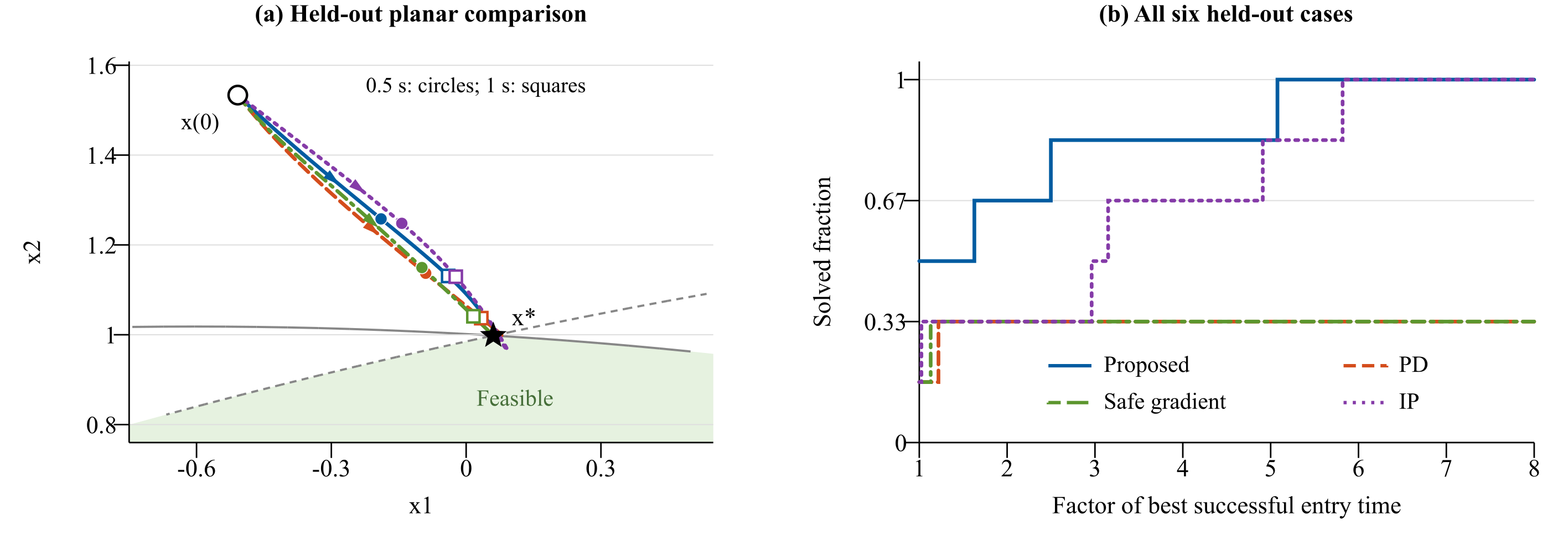}
 \caption{Test results at $U=4$. Left: trajectories on the first planar
 test case, with the feasible region shaded and the optimizer marked
 by a star. Right: entry-time profile over all six cases.}
 \label{fig:comparison}
\end{figure*}

Each run is evaluated on $[0,20]$. The common violation allowance
$\bar\rho_{2,i}(t):=.02+(w_i-.02)(1-t/4)_+^2$, with
$w_i:=\max\{1,.5+1.2g_i^+(x_0)\}$, reaches $.02$ at $t=4$ and remains
there. The entry time $t_{.02}$ is the first recorded time from which
$\norm{x-x^\star}<.02$ holds at all remaining samples. A run succeeds
if it reaches $t=20$, enters the \(0.02\) neighborhood of the optimizer and stays there at all subsequent recorded samples, and satisfies
$g_i^+(x(t))/\bar\rho_{2,i}(t)\leq1+10^{-6}$ and
$\norm{\dot x(t)}\leq U$ at every sample, for the speed cap
$U\in\{2,4,8\}$. Primal effort
$E_x:=\int_0^{20}\norm{\dot x(t)}^2dt$ is estimated by quadrature.

Projected PD~\cite{feijer2010} uses $\dot x=-a_xr(x,\lambda)$ and
$\dot\lambda_i=a_\lambda[g_i(x)]_{\lambda_i}^+$, with
$a_x,a_\lambda>0$. Here $[g_i(x)]_{\lambda_i}^+$ equals $g_i(x)$ for
$\lambda_i>0$ and $g_i^+(x)$ for $\lambda_i=0$.
Safe gradient~\cite{allibhoy2024} minimizes
$\norm{v+\nabla f(x)}^2/2$ over $v\in\R^n$ subject to
$J_g(x)v\leq-a_gg(x)$, and sets $\dot x=a_vv$, with $a_g,a_v>0$.
Prediction--correction IP~\cite{fazlyab2018} solves a Newton system
with correction gain $a_{\rm IP}>0$ for the barrier objective $ \Phi(t,x):=f(x)-c(t)^{-1}\sum_i\log(\sigma(t)-g_i(x)),$ where $\sigma(t):=\sigma_0e^{-\eta_st}$, $c(t):=c_0e^{\eta_ct}$,
$\sigma_0>\max\{0,\max_i g_i(x_0)\}$, and $\eta_s,c_0,\eta_c>0$.
The proposed update is explicit and first order.

The tuning budget is one default and 31 deterministic radical-inverse
configurations per method, evaluated on the same development
problems~\cite{beiranvand2017}. Both PD~\cite{feijer2010} gains lie in
$[.2,8]$. For safe gradient~\cite{allibhoy2024},
$(a_v,a_g)\in[.25,8]\times[.2,4]$. For IP~\cite{fazlyab2018}, the ranges
of $(a_{\rm IP},\kappa,c_0,\eta_s,\eta_c)$ are $[.25,6]$, $[1.1,4]$,
$[1,40]$, $[.4,2.5]$, and $[.08,.5]$, respectively, with
$\sigma_0:=\kappa\max\{1,\max_i g_i(x_0)\}$.
The proposed method uses $k_1,k_\lambda\in[.2,4]$,
$\rho_{2,i}=\bar\rho_{2,i}$, and
\begin{equation}\label{eq:comparison-profiles}
 \begin{aligned}
 \rho_{1,j}(t)&:=\kappa_r(.25+|[\nabla f(x_0)]_j|)
                      \big(h_r+(1-h_r)e^{-b_rt}\big),\\
 \rho_S(t)&:=\kappa_S\Big(1+\sum_i s_i(0)^2\Big)
                      \big(h_S+(1-h_S)e^{-b_St}\big).
 \end{aligned}
\end{equation}
The ranges are $(\kappa_r,h_r,b_r)\in[1.3,3]\times[.02,.2]\times[.2,1.2]$
and $(\kappa_S,h_S,b_S)\in[.3,3]\times[.001,.1]\times[.2,1]$.
No Slater or curvature data are used. For each method and speed cap,
the configuration with the most development successes is selected.
Ties favor lower mean entry time (assigning $40$ to a failure), then
lower maximum sampled speed, then configuration ID.

\begin{table}[t]
 \caption{Test successes out of six cases. At $U=4$, $t_{.02}$ and
 $E_x$ are means over successful runs. A dagger marks a configuration
 with only two development successes.}
 \label{tab:comparison}
 \centering\footnotesize
 \setlength{\tabcolsep}{3pt}
 \begin{tabular}{lrrrrr}
 \hline
 Method & $U=2$ & $U=4$ & $U=8$ & $t_{.02}$ & $E_x$\\
 \hline
 Proposed & $4/6$ & $6/6$ & $6/6$ & $6.050$ & $.562$\\
 PD~\cite{feijer2010} & $2/6^\dagger$ & $2/6$ & $2/6$ & $1.504$ & $.848$\\
 Safe gradient~\cite{allibhoy2024} & $3/6^\dagger$ & $2/6$ & $4/6$ & $1.432$ & $.869$\\
 IP~\cite{fazlyab2018} & $3/6$ & $6/6$ & $6/6$ & $6.835$ & $.678$\\
 \hline
 \end{tabular}
\end{table}

At $U=4$, the proposed method and IP~\cite{fazlyab2018} succeed on all
six test cases, whereas PD~\cite{feijer2010} and safe
gradient~\cite{allibhoy2024} succeed on two each
(Table~\ref{tab:comparison}). Over the six shared successes, the proposed
method has lower mean entry time and primal effort than
IP~\cite{fazlyab2018}. The two-case means do not support an overall
speed ranking. Figure~\ref{fig:comparison} shows the trajectories on one
test case and the entry-time profile~\cite{dolan2002}, which assigns
infinite cost to failures.

\section{Conclusion}\label{sec:conclusion}
The proposed dynamics enforce prescribed transient bounds and converge
asymptotically to the constrained optimizer. With known objective
curvature and strict-feasibility data, a prescribed accuracy can also be
guaranteed from a chosen deadline onward. Potential future work includes applications to robotic control allocation and energy dispatch, where temporary
violations of soft operating limits can be bounded while new actuator
references or dispatch setpoints are computed.

\bibliographystyle{IEEEtran}
\bibliography{reference}
\end{document}